\input amssymb.sty

\documentclass{article}

\usepackage{url}

\newtheorem{Proof}{Proof.}

\newtheorem{theorem}{Theorem}

\newtheorem{ax}{A}

\newenvironment{proof}{\begin{Proof}\rm}{\hfill $\Box$ \end{Proof}}

\begin{document}

\title{The sum of irreducible fractions with
consecutive denominators is never an integer in a very weak arithmetic}

\author{Victor Pambuccian}

\date{}

\maketitle

Most problem solvers have encountered at some stage of their lives
the problem asking for a proof that the sum

\[ 1+\frac{1}{2}+\ldots + \frac{1}{n}\]
for $n\geq 2$ is never an integer.  The proof one usually finds
offered for this fact is based on Chebyshev's theorem (Bertrand's
postulate). If one asks for a proof that, more generally, the sum

\[ \frac{1}{n}+\frac{1}{n+1}+\ldots + \frac{1}{n+k}\]
with $k\geq 1$ can never be an integer, then the proof based on Chebyshev's
theorem needs to be amended. One first notes that, if $k<n$, then the above
sum must be less than 1, and thus cannot be an integer, and if $k\geq n$,
then one applies the same proof based on Chebyshev's theorem (this
fact seems to have been overlooked in \cite{obl}, where the
author  wants to use Chebyshev's Theorem, but
finds it necesary to make the proof dependent on
another deep result, the Sylvester-Schur Theorem, as well).
However, K\"ursch\'ak \cite{kuer} (see also \cite{ps}) found a
much simpler proof, which relies on the very simple observation
that among any number ($\geq 2$)
of consecutive positive integers there is precisely one, which is divisible
by
the highest power of 2 from among all the given numbers. Aside from its
didactical
use, one may wonder whether K\"ursch\'ak's proof is not in a very formal
way much simpler, i.\ e.\ whether it does not require simpler methods of
proof
in the sense of formal logic.

When formalized, arithmetic is usually presented as Peano
Arithmetic, which contains an induction axiom schema, stating,
loosely speaking, that any set that can be defined by an elementary
formula in the language of arithmetic (i.e. in terms of some
undefined operation and predicate symbols, such as $+$, $\cdot$,
$1$, $0$, $<$), which contains 1, and which contains $n+1$ whenever
it contains $n$, is the set of all numbers.  Several weak
arithmetics have been studied, in which the types of elementary
formulas allowed in the definitions of the sets used in induction
are restricted by certain syntactic constraints (see \cite{daq}),
and one might think that K\"ursch\'ak's proof would make it in a
weaker formal arithmetic than the one dependent on Chebyshev's
theorem. It turns out that, in fact, no amount of induction is
needed at all!

To see this, let's first generalize the problem further, along the
lines of the generalization in \cite{obl1}, so that there can be no
proof based on Chebyshev's Theorem.

\begin{theorem}

The sum

\begin{equation} \label{k1}
\frac{m_0}{n}+\frac{m_1}{n+1}+\ldots + \frac{m_k}{n+k}
\end{equation}

with $(m_i,n+i)=1$, $m_i<n+i$, and $k\geq 1$ is never an integer.
\end{theorem}

\begin{proof} (K\"ursch\'ak \cite{kuer}). Let $a = \max \{ \alpha :
2^{\alpha} | (n+i) \mbox{  for some
$0\leq i \leq k$ } \}$. Then $2^a$ divides exactly one of the numbers
$n$, $n+1$, \ldots, $n+k$. Let $l= \mbox{lcm } (n, n+1, \ldots, n+k)$.
Suppose the sum in (\ref{k1}) is an integer
$b$. Multiplying both (\ref{k1}) and $b$ by $l$, we obtain on the one hand
an odd number, and on the other an even number,
which have to be equal.
\end{proof}

Moreover, to make it a theorem of arithmetic, we will do away with
the fractions appearing in it, and state it, for all positive $k\in
\mathbb{N}$, as $\varphi_k$, the following statement (where we
denote by $\overline{u}$ the term $((\ldots ((1+1)+1)+\ldots)+1)$,
in which there are $u$ many $1$'s; the terms $\overline{u}$ will be
referred to as {\em numerals})
\begin{eqnarray*}
& & (\forall n)(\forall m_0)\ldots (\forall m_k)(\forall p)\,
\bigvee_{i=0}^k ((\forall a)(\forall b)\, m_ia \neq (n+\overline{i})b +1)\vee \bigvee_{i=0}^k n+\overline{i}<m_i\\
& & \vee \sum_{i=0}^k (m_i\prod_{0\leq j\leq k, j\neq i}
(n+\overline{j})) \neq p \prod_{j=0}^k (n+\overline{j}).
\end{eqnarray*}

The arithmetic we will show it holds in is PA$^-$, which is
expressed in a language containing as  undefined operation and
predicate symbols only $+$, $\cdot$, $1$, $0$, and $<$, and whose
axioms A1-A15 were presented in \cite[pp.\ 16-18]{kay}. We will
repeat them here for the reader's convenience, and we will omit the
universal quantifiers for all universal axioms.

\begin{ax}
$(x+y)+z = x+(y+z)$
\end{ax}

  \begin{ax}
$x+y =  y+x$
\end{ax}

  \begin{ax}
$(x\cdot y)\cdot z = x\cdot (y\cdot z)$
\end{ax}

  \begin{ax}
$x\cdot y = y\cdot x$
\end{ax}

  \begin{ax}
$x\cdot (y+z) = x\cdot y+x\cdot z$
\end{ax}

  \begin{ax}
$x+0=x \wedge x\cdot 0 = 0$
\end{ax}

  \begin{ax}
$x\cdot 1 = x$
\end{ax}

  \begin{ax}
$(x<y \wedge y<z) \rightarrow x<z$
\end{ax}

  \begin{ax}
$\neg x<x$
\end{ax}

  \begin{ax}
$x<y \vee x=y \vee y<x$
\end{ax}

  \begin{ax}
$x<y \rightarrow x+z < y+z$
\end{ax}

  \begin{ax}
$(0<z \wedge x<y) \rightarrow x\cdot z < y\cdot z$
\end{ax}

  \begin{ax}
$(\forall x)(\forall y)(\exists z)\, x<y \rightarrow x+z =y$
\end{ax}

  \begin{ax}
$0<1 \wedge (x>0 \rightarrow (x>1 \vee x=1))$
\end{ax}

  \begin{ax}
$x>0 \vee x=0$
\end{ax}

What is missing from $PA^-$, and makes it so weak (indeed, the
positive cone of every discretely ordered ring is a model of
$PA^-$), is the absence of any form of induction.

The proof that  $\varphi_k$ holds in $PA^-$ will be carried out in
an arbitrary model $\mathfrak{M}$ of $PA^-$. The idea of proof will
be to show that all variables that appear in $\varphi_k$ must be
numerals. An essential ingredient of the proof is the following
fact, which holds in $PA^-$ (see \cite[Lemma 2.7, p.\ 22]{kay}), for
all positive $k\in\mathbb{N}$
\[x<\overline{k} \rightarrow x=0\vee x=1\vee \ldots \vee x=\overline{k-1},\]
and which allows us to deduce that any element that is bounded from
above by a numeral must be a numeral.

Suppose that, for some positive $k\in \mathbb{N}$, $\varphi_k$ does
not hold in $\mathfrak{M}$. Then, for all $i=0, \ldots, k$, there
are $m_i, p, a_i$  and $b_i$ with $m_ia_i=(n+i)b_i + 1$ and such that
\begin{equation} \label{1}
\sum_{i=0}^k (m_i\prod_{0\leq j\leq k, j\neq i} (n+\overline{j})) = p
\prod_{j=0}^k
(n+\overline{j}).\end{equation}

This can be rewritten, by leaving only the first term of the sum on
the left-hand side, and sending all others to the right-hand side
with changed sign,   as $m_0(n+\overline{1})\ldots (n+\overline{k})=
nq$, where by $q$ we have denoted $p\prod_{j=1}^k (n+\overline{j}) -
(\sum_{i=1}^k m_i\prod_{0\leq j\leq k, j\neq i} (n+\overline{j}))$.
The product $(n+\overline{1})\ldots (n+\overline{k})$ can also be
written as a polynomial in $n$, whose free term is $\overline{k!}$,
i.\ e.\ as $nr + \overline{k!}$, thus $m_0(nr + \overline{k!})=nq$.
Given that there are  $a_0$ and $b_0$ such that  $m_0a_0=nb_0
+1$, if we multiply both sides of the equality  $m_0(nr + \overline{k}!)=nq$
by $a_0$ we obtain $(nb_0+1)\overline{k!} = n(a_0q-a_0m_0r)$,
thus $\overline{k!} = n(a_0q-a_0m_0r-b_0\overline{k!})$.
We know that $\mathfrak{M}$ must contain a
copy of $\mathbb{N}$, and it may contain other elements as well,
called {\em nonstandard} numbers.  Could $n$ be in $\mathfrak{M}$
but not of the form $\overline{m}$ for some $m\in \mathbb{N}$? If it
were such an element of $\mathfrak{M}$, then it would be greater
than all $\overline{m}$ with $m\in \mathbb{N}$, and thus so would
$n(a_0q-a_0m_0r-b_0\overline{k!})$, unless $a_0q-a_0m_0r-b_0\overline{k!}=0$,
which cannot be the case, as $\overline{k!}$
is not zero. However, $n(a_0q-a_0m_0r-b_0\overline{k!})$ cannot be greater than all
$\overline{m}$ with $m\in \mathbb{N}$, for it is equal to such a
number, namely to $\overline{k!}$. Thus $n$ must be an
$\overline{m}$ for some $m\in \mathbb{N}$.  This means that in
(\ref{1}) all variables are numerals, i.\ e.\   $\overline{i}$'s for
some $i\in \mathbb{N}$. However, we know, from K\"ursch\'ak's proof,
that such an equation cannot exist, so, for all $k\in \mathbb{N}$,
$\varphi_k$ holds in PA$^-$.

Another generalization of the original problem, proved by T. Nagell
in \cite{nag}, states that the sum
\[\frac{1}{m}+\frac{1}{m+n}+\frac{1}{m+2n}+\cdots+\frac{1}{m+kn}\]
is never an integer if $n, m, k$ are positive integers. The proof is
rather involved and uses both a K\"ursch\'ak-style argument and
Chebyshev's theorem.  This statement turns out to be, with
$\overline{k}$ instead of $k$, valid in $PA^{-}$ as well. To see
this,  let, for all positive $k\in \mathbb{N}$,  $\nu_k$ stand for

\begin{equation} \label{nagell}
(\forall m)(\forall n)(\forall p)\, m>0 \wedge n>0 \rightarrow
\sum_{i=0}^{k}\prod_{0\leq j\leq k, j\neq i} (m+\overline{j}n) \neq
p\prod_{0\leq j\leq k} (m+\overline{j}n),
\end{equation}
and let $\mathfrak{M}$ be again a model of $PA^{-}$. Notice that, if
$m> \overline{k}$, then, then $\neg\nu_k$ cannot hold for any $n$
and $p$. To see this, suppose that, for some $n$ and $p$, we have
equality in (\ref{nagell}). Given that $(m+n)\ldots
(m+\overline{k}n)$ is the largest of all the summands on the
left-hand side, and there are $\overline{k}$ summands, the sum on
the left-hand side is $\leq \overline{k}(m+n)\ldots
(m+\overline{k}n)$, and thus $< m(m+n)\ldots (m+\overline{k}n)$,
thus equality cannot hold in (\ref{nagell}). Thus $m\leq
\overline{k}$, and thus (see \cite[p.\ 20]{kay}) $m$ must be
standard, i.\ e.\ it must be $\overline{u}$ for some $0<u\leq k$. It
remains to be shown that $n$ must be standard as well. To see this,
suppose again that, for some $n>0$ and $p$, we have equality in
(\ref{nagell}). Notice that, since $m(m+2n)\ldots (m+\overline{k}n)$
is the largest product among all $\prod_{0\leq j\leq k, j\neq i}
(m+\overline{j}n)$, for $i=1,2, \ldots k$, the sum on the left hand
side of our equality is $\leq (m+n)\ldots (m+\overline{k}n) +
\overline{k}m(m+2n)\ldots (m+\overline{k}n)$ (with equality if and
only if $k=1$). Thus $pm(m+n)\ldots (m+\overline{k}n) \leq (m+n +
\overline{k}m)(m+2n)\ldots (m+\overline{k}n)$, which implies
$pm(m+n)\leq m+n + \overline{k}m$. If $n$ were nonstandard, then
this inequality were possible only if $pm=1$, i.\ e.\ if $p=m=1$,
which is not possible, for in that case the first summand on the
left hand side of (\ref{nagell}) is equal to the right hand side,
thus the left hand side must be larger than the right hand side, so
equality could not have taken place in (\ref{nagell}). Now that
$m,n,k$ have all been shown to be standard, Nagell's proof implies
the truth of our statement, which thus holds in $PA^{-}$.

By G\"odel's completeness theorem, there must exist  syntactic
proofs that $\varphi_k$ and $\nu_k$ hold in PA$^-$, i.\ e.\  formal
derivation of $\varphi_k$ and $\nu_k$ from the axioms of  PA$^-$.
Such a formal proof for $\varphi_k$ cannot use the idea behind
K\"ursch\'ak's proof, for it is not even true in  PA$^-$ that among
$n$, $n+1$, $\ldots$ $n+k$, there is a multiple of 2 (see \cite[p.\
18]{kay}). Thus there must  exist even simpler proofs for both
$\varphi_k$ and $\nu_k$ and they are worth finding. Such proofs
would reveal the real reasons why these results hold, and the reason
must be of an algebraic nature, for there is no traditional number
theory to be found in $PA^-$.

\bigskip

\noindent\textit{Department of Integrative Studies, Arizona State
University - West Campus, Phoenix, AZ 85069-7100, U.S.A\\
pamb@math.west.asu.edu}


\begin{thebibliography}{17}


\bibitem{daq} P. D'Aquino, Weak fragments of Peano Arithmetic,
in: \textit{The Notre Dame Lectures}, (ed. P. Cholak), Lecture Notes
in Logic, 18, pp. 149--185, A. K. Peters, 2005.


\bibitem{kay} R. Kaye, \textit{Models of Peano Arithmetic}, Oxford
University Press,
Oxford, 1991.


\bibitem{kuer} J.  K\"ursch\'ak, On the harmonic series (Hungarian),
Mat.\ Phys.\ Lapok 27 (1918), 299--300.


\bibitem{nag} T. Nagell, Eine Eigenschaft gewisser Summen,
Norsk Mat. Forenings Skrifter 13 (1923), 10--15.

\bibitem{obl1} R. Oblath, Some number-theoretical theorems (Hungarian),
Math. \'es phys. lapok 27 (1918), 91--94.


\bibitem{obl} R. Oblath, \"Uber einen arithmetischen Satz von K\"ursch\'ak,
Comment.\ Math.\ Helv.\  8  (1935),   186--187.


\bibitem{ps} G. P\'olya, G. Szeg\"o, \textit{Aufgaben und Lehrs\"atze der
Analysis},
Springer-Verlag, Berlin, 1971, vol. II, VIII. Abschnitt, Aufgabe
251, S.\ 159.







\end{thebibliography}
\end{document}